\documentclass[a4paper,10pt]{article}

\usepackage[a4paper,left=3cm,right=3cm,top=2.5cm,bottom=2.5cm]{geometry}
\usepackage{hyperref}
\usepackage{subfig}
\usepackage{pgfplots}
\usepackage{pgfplotstable}
\usepackage{mathtools}
\usepackage{multicol}
\usepackage{comment}
\usepackage{booktabs}
\pgfplotsset{compat=1.5}
\usepackage{amssymb}
\usepackage{url}
\usepackage{bm}

\usepackage{pdfsync}
\usepackage{float}
\usepackage{tabularx}
\usepackage{enumerate}
\usepackage{array}
\usepackage{xspace}
\usepackage{siunitx}
\usepackage{authblk}

\usepackage{amsthm}
\usepackage{amsmath}
\usepackage{stmaryrd}
\usepackage[ruled,vlined,linesnumbered]{algorithm2e}
\usepackage{graphicx}
\usepackage{booktabs}
\usepackage{cleveref}

\usepackage{overpic}

\newtheorem*{remark}{Remark}
\DeclareMathOperator*{\argmin}{\arg\!\min}

\usepackage[inline, shortlabels]{enumitem}
\setlist{nolistsep}

\newcommand{\R}{\mathbb R}

\newcommand{\bI}{\mathbf I}

\newcommand{\btau}{\boldsymbol{\tau}}

\newcommand{\bd}{\mathbf d}

\newcommand{\bn}{\mathbf n}

\newcommand{\bu}{\boldsymbol{u}}

\newcommand{\bx}{\mathbf x}

\newcommand{\B}{\mathcal B}

\newcommand{\J}{\mathcal J}

\newcommand{\bff}{{\boldsymbol{f}}}
\newcommand{\bur}{{\boldsymbol{u}}_r}

\newcommand{\bvr}{{\boldsymbol{v}}_r}
\newcommand{\bphi}{\boldsymbol{\varphi}}
\newcommand{\bX}{\mathbf X}
\newcommand{\bXr}{{\mathbf X}^r}

\begin{document}

\title{Digital Twins in Coronary Artery Disease: \\ A Mathematical Roadmap}

\author[1]{Alessandro Veneziani\footnote{avenez2@emory.edu}}
\author[2]{Annalisa Quaini\footnote{aquaini@central.uh.edu}}
\author[1]{Marco~Tezzele\footnote{marco.tezzele@emory.edu}}
\author[3]{Omer San\footnote{osan@utk.edu}}
\author[4]{Traian Iliescu\footnote{iliescu@vt.edu}}

\affil[1]{Emory University, 400 Dowman Drive NE, Atlanta, GA 30322, United States}
\affil[2]{University of Houston, 3551 Cullen Blvd, Houston, TX 77204, United States}
\affil[3]{University of Tennessee, 1512 Middle Drive, Knoxville, TN 37996, United States}
\affil[4]{Virginia Tech, 225 Stanger Street, Blacksburg, VA 24061, United States}

\maketitle

\vspace{-.2cm}

\begin{abstract}
The combination of data and models, enhanced by AI methodologies, leads to the paradigm called \textit{Digital Twins}. This concept is expected to bring unprecedented support to personalized medicine. The combination of mathematical and numerical models with diagnostic devices that provide patient-specific knowledge in a bidirectional framework can be a formidable decision support for clinicians. In this paper, we consider some mathematical aspects of constructing a Digital Twin to prevent and treat Coronary Artery Disease. The keywords for the bidirectional communication between twins 
in our system 
are (i) Data Assimilation and (ii) Probabilistic Graphic Models. In particular, a quantity of paramount interest in the evaluation and prognosis of Coronary Artery Disease is the Wall Shear Stress, i.e., the tangential component of normal stress on the arterial wall. By considering steps for the personalization and the synthesis of Wall Shear Stress estimation, we propose a mathematical roadmap for constructing a Digital Twin system that could help prevent infarcts, one of the most lethal diseases in the world.
\end{abstract}

\tableofcontents

\section{Digital-Twins for Cardiovascular Health: Landscape Scenario}
\label{sec:intro}

With the advent of data-driven models and AI methodologies, the construction of virtual models of physical and engineering systems,
previously based on computational mechanics, has ascended to unprecedented levels of complexity. Not only data-driven modeling enables the use of virtual descriptions of complex systems in a wider range of fields (biology, sociology, urban development, etc.), but the interaction between the physical and 
the virtual systems is brought to a much more sophisticated paradigm that takes the name of ``digital twins" (DTs). 
The following definition 
is provided in 
\cite{nationalacademy}:

\begin{quotation}
    A digital twin is a set of virtual information constructs that mimics the structure, 
context, and behavior of a natural, engineered, or social system (or system of systems), is dynamically updated with data from its physical twin, has a 
predictive capability, and informs decisions that realize value. The bidirectional 
interaction between the virtual and the physical is central to the digital twin.
\end{quotation}

Beyond the predictive nature of virtual models, 
the (bidirectional) integration of data and conceptualizations
opens up a wealth of perspectives. The virtual twin
is not only a tool for understanding the physical system and designing
a better one, but also a ``continuous-in-time'' decision-making support. One of the most promising and impactful fields where DTs are expected to make a difference is medicine, as demonstrated by the abundant literature in the field \cite{dziopa2024digital,sel2024building,sel2025survey}. 
In particular, cardiovascular diseases are of interest for at least two different reasons: (i) they are still the most lethal diseases in the Western world \cite{cdcFastStats}; (ii) mathematical modeling of the cardiovascular system has an exceptional level of maturity and depth (see, e.g., \cite{formaggia2010cardiovascular,taylor1998finite,quarteroni2000computational,quarteroni2017cardiovascular}).

Although the perspectives of DTs for understanding and curing cardiovascular diseases are enormous \cite{bruynseels2018digital,corral2020digital,coorey2022health,sel2024building,dziopa2024digital},
many aspects remain challenging for several reasons. 
In particular,
as noted in \cite{sel2024building}, mathematical models are well developed and advanced, yet they lack a rigorous personalization.
To be more precise, 
the cardiovascular apparatus features a hierarchy of components,
{which 
correspond to different physical scales. 
This hierarchy goes from the molecular level (critical for understanding the
electrophysiology of the heart), to the cellular level, the tissue level, and the organ level, up to the systemic level.
The interplay of these different scales has been a challenge for applied mathematicians for decades, 
and led to the introduction of specific numerical methods, see, e.g., \cite{formaggia1999multiscale,vignon2010outflow,quarteroni2016geometric,perdikaris2016multiscale}. On the other hand, 
the data collection and usage pose significant challenges,
both from the technological and the methodological point of view. In cardiovascular medicine, there are different types of data collected at different 
{frequencies} and, thus, the ``dynamical update'' advocated by the definition of DT has different implementations.
We may have data from (i) \textit{wearable technologies}, e.g., smart rings or smart watches, with measurements almost continuous in time, (ii) \textit{ portable technologies}, e.g., pressure measurement devices, with daily (or more frequent) measurements, like home 
sphygmomanometers; (iii) \textit{ noninvasive devices}, e.g., ultrasound (US), that can be collected rather frequently in principle, and finally, (iv) \textit{invasive systems}, e.g., computed tomography (CT) and magnetic resonance imaging (MRI), which can be collected only sporadically. 
In general, the quality of the data is inversely proportional to their frequency.
The heterogeneity of the data and their frequency of collection
{are} critical aspects to address, as well as the intrinsic ``multiscale-in-time'' nature of many cardiovascular diseases, where the time scale of a heart beat (order of seconds) impacts the development of the
pathology over weeks or months.

All these aspects challenge applied mathematicians and urge us to extend the concept of ``mathematics for clinics'' to unexplored territories.
Currently, some of the most important keywords of this process are: (i) personalization and (ii) synthesis. 
They concern the bidirectionality of the DT system, being central both at the input (how do the data make a general model patient-specific?)
and at the output (how do the results of the model propagate to the patient?) of the virtual twin.
The concept of ``personalization'' may sound more in the realm of engineering than of mathematics, as it is antithetical to the
abstraction of a mathematical model, and of providing proofs of concept
at a general level. 
However, in our opinion, these steps (personalization and synthesis) present a formidable number of challenges for the combination of 
deterministic and stochastic components, {and} deductive and inductive reasoning. 
Only the best mathematical tools can provide trustworthy answers to such challenges.

In this paper, we focus on the development of a DT for Coronary Artery Disease (CAD), depicted in Fig.~\ref{fig:cad_scheme}.
The reason for this specific interest is twofold.
First, CADs are 
the most common heart disease. Heart diseases are responsible of more than $700{,}000$
exitus per year in the USA \cite{cdcFastStats}. 
Second, the medical community agrees on the role of a specific hemodynamic quantity, the so-called Wall Shear Stress (WSS), as a major actor in the development of CAD. Although the WSS 
is critical, it cannot be directly measured. 
To address this challenge, in this paper we  
focus on the mathematical personalization of WSS estimation as a combination of data and models, and on the  consequent synthesis to inform the decision-maker. 
We aim at identifying a realistic roadmap for designing a first DT that could have a clinical impact on the prediction
and prevention of severe adverse events related to CAD, e.g., infarction. 
This paper focuses on the robust mathematical foundation needed for such a DT to become a reality in clinical practice.

\begin{figure}[htb!]
\centering
\includegraphics[width=\textwidth]{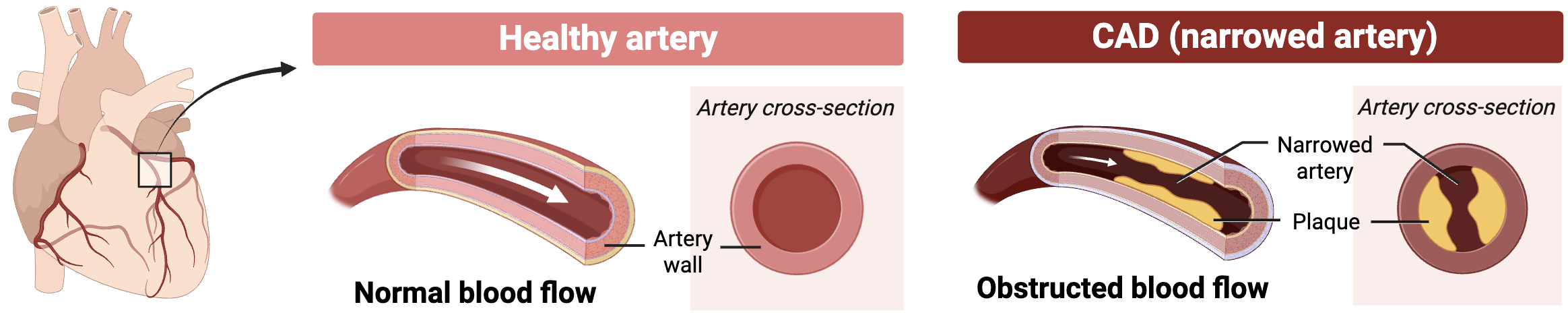}
\caption{Comparison between
a healthy coronary and a coronary 
affected by CAD. Image created in \url{bioRender.com}.}
\label{fig:cad_scheme}
\end{figure}

The paper is organized as follows. We start by giving a short introduction to CAD at the mathematical and 
clinical levels. Sec.~\ref{sec:dt} describes the general organization of a DT system.
Sec.~\ref{sec:VDA} considers possible techniques for the personalization of the WSS estimation, based on velocity measures that can be retrieved from noninvasive devices like US.
The keyword in this stage is Data Assimilation (DA), and we consider different ways to enable it.
Then, Sec.~\ref{sec:PGM} discusses methods of synthesis to organize the data into decision-making support focusing, in particular, on Probabilistic Graphic Models (PGMs).
Finally, we give some assessments of the perspectives of DTs for CAD and beyond in Sec.~\ref{sec:road}.

\section{Basics on Coronary Artery Disease}\label{sec:back}

In CAD, abnormal calcifications and lipid deposits create occlusions in the coronary arteries (see Fig.~\ref{fig:cad_scheme})
that can eventually lead to myocardial infarction. This section provides a brief presentation of CAD from a mathematical perspective in Sec.~\ref{subsec:math_back} and from a clinical point of view in Sec.~\ref{subsec:clinical_back}.

\subsection{The Mathematical Framework}
\label{subsec:math_back}
CAD stems from a complex interplay between blood and the vascular tissues. A comprehensive mathematical model for this pathology should include the different physical components (the fluid dynamics, the structural dynamics, and the biological processes)
and cover the different relevant time scales. 
In this design of a first decision-making support system,  
we focus specifically on the hemodynamics. Consequently, we cover here only the aspects related to the fluid.

The blood flow in the coronary arteries is modeled by the incompressible Navier-Stokes equations 
(NSE) for a Newtonian fluid in laminar regime\footnote{In the presence of strong occlusions, 
the flow regime can be highly disturbed. In general, however, in coronary hemodynamics an appropriate mesh
and, possibly, classical stabilization techniques are enough for the numerical simulations.}. Even 
though this includes simplifications (for instance, the viscosity is generally non constant), this is a commonly accepted model.
We denote by $\Omega$ the domain occupied by the blood, which is a bounded open set in $\R^3$ with 
Lipschitz boundary $\Gamma = \partial \Omega$. We denote by $\bu$ the velocity of the blood, 
by $p$ its pressure,
by $\bff$ the body force per unit mass, 
by $\mu >0$ the dynamic viscosity of the blood, 
and by $\rho >0$ its density. The stress tensor $\btau$ is defined as
\begin{equation}
    \btau = -p \bI + \mu \left( \nabla \bu + (\nabla \bu)^T \right), 
\end{equation}
where $\bI$ is the identity tensor. Then,  the NSE read as follows: find $\bu$ and $p$ such that
\begin{eqnarray}
    \rho \left( \frac{\partial \bu}{\partial t} + (\bu \cdot \nabla)\bu \right) - \nabla \cdot \btau &= \bff \quad \text{in } \Omega \times (0,T), \label{eq:pbc_nse} \\
    \nabla \cdot \bu &= 0 \quad \text{in } \Omega \times (0,T). \label{eq:pbc_nse-2}
\end{eqnarray}

These equations should be complemented with suitable initial and boundary conditions, 
depending on the specific problem setting.
For the initial conditions, generally a number of heart beats are run to reach a periodic regime, and to lose dependence on an arbitrary initial condition (which can be computed as the solution of a steady Stokes problem
in $\Omega$). Regarding the boundary conditions,  it is reasonable
to assume that blood flow is driven by a pressure drop across the coronary artery.
So, we prescribe no-slip boundary conditions on the arterial wall, i.e.,
\begin{equation}
    \bu = \mathbf{0} \quad \text{on } \Gamma_{\text{wall}} \times (0,T),
\end{equation}
where $\Gamma_{\text{wall}}$ denotes the portion of the boundary corresponding to the vessel wall. 
At the inlet and outlet sections, we prescribe stress boundary conditions including the pressure,
\begin{eqnarray} 
    \btau_{{\text{in}}} \cdot \bn &= -p_{\text{in}} \bn \quad &\text{on } \Gamma_{\text{in}} \times (0,T), \label{eq:pbc} \\
    \btau_{{\text{out}}} \cdot \bn &= -p_{\text{out}} \bn \quad &\text{on } \Gamma_{\text{out}} \times (0,T), \label{eq:pbc-2}
\end{eqnarray}
where $p_{\text{in}}$ and $p_{\text{out}}$
 are the prescribed (possibly time-dependent) pressures at the inlet and outlet, respectively, 
 $\bn$ is the outward normal vector to the boundary, 
 and $\Gamma_{\text{in}}$ and $\Gamma_{\text{out}}$ are the inflow and outflow portions of the boundary, respectively.
 However, we note that  in practice the pressure can be measured only invasively and these data are generally
 not available. For this reason, these conditions can be replaced by conditions on the velocity. When not available, the velocity can be inferred by data on the flow rate, 
 or other geometrical considerations, with reasonable assumptions on the profile of the velocity itself. We will, however, refer to boundary conditions~\eqref{eq:pbc}, as they are functional in our DA procedure, as we illustrate later on.

As mentioned in the introduction, a quantity of great interest in the assessment of CAD is the WSS, defined as the tangential component of the stress tensor $\btau$ on the arterial wall:
\begin{equation}
    \text{WSS} \equiv \btau \cdot \bn - (\bn\cdot \btau \cdot \bn)\bn \quad \text{on } \Gamma_{\text{wall}} \times (0,T) \,.    
\end{equation}

Many quantities derived from the WSS are used for predictions value.
Among the others, we mention the Oscillatory Shear Index~\cite{ku1985pulsatile},
the Relative Residence Time~\cite{murray1926physiological}, and the Time-Averaged WSS (TAWSS)~\cite{ku1985pulsatile}.
The latter is simply the time-averaged value of the WSS over a cardiac cycle:
\begin{equation}
    \text{TAWSS} = \frac{1}{T} \int_0^T | \text{WSS} (t)|\, dt \,,
\end{equation}
where $T$ is the period of the cardiac cycle.

The WSS (its magnitude or the TAWSS) can be  estimated by 
post-processing measured velocities. In \cite{martin2024pseudo}, 
a pseudo-spectral approach that filters high-frequency noise in ultrasound velocity measures 
and provides more accurate estimates
than traditional Finite Difference formulas is proposed. This is a purely data-driven approach.
Another, model-driven, approach is based on computational fluid dynamics (CFD), where the velocity field is computed by solving numerically the NSE. This requires knowledge of data on the boundaries. Unfortunately, most of the data that can be measured, in practice, are not on the boundary, so  cannot be trivially used to inform the boundary conditions. 
For the TAWSS, a commonly accepted approach is to solve the steady Navier-Stokes equations
and then to compute directly the TAWSS from the steady velocity field $\bu$.
Recent quantitative studies have shown that this is an acceptable surrogate for the TAWSS
\cite{shah2025impact,nannini2024automated}. 

This second CFD-based approach suffers from the lack of boundary data and their uncertainty. 
In our DT framework, we propose a more complex (and hopefully more accurate) approach,
obtained by merging the previous two strategies (Sec. \ref{sec:VDA}).

\subsection{The Clinical Framework}
\label{subsec:clinical_back}

Chest pain is the symptom of CAD that usually leads a patient to the Emergency Room, where the severity of the disease can be
assessed by measuring hemodynamic indices \cite{fearon2016fractional,samady2018interventional,clevelandclinicFractionalFlow}. 
The most popular index is called (myocardial) \textit{Fractional Flow Reserve} (FFR) 
\cite{pijls1993experimental,de1994coronary,pijls1995fractional}. The FFR estimates 
the functional significance of a stenosis as the quotient between the maximum achievable blood flow in the presence of a stenosis 
over the one in physiological conditions; referring to the Bernoulli ``vis-viva'' equation \cite{formaggia2010cardiovascular}, this quantity is expressed in terms of pressures as:  
\begin{equation}\label{eq:ffr}
    \text{FFR} = \frac{p_{\text{distal}}-p_{\text{venous}}}{p_{\text{arterial}}-p_{\text{venous}}} =
    1 - \frac{\delta P}{p_{\text{arterial}}-p_{\text{venous}}} \,,
\end{equation}
where $p_{\text{distal}}$ is the pressure downstream of the stenosis, $p_{\text{venous}}$ is the pressure in the right ventricle, 
$p_{\text{arterial}}$ is the pressure upstream of the stenosis, 
and $\delta P$ is the trans-stenotic pressure
(see Fig.~\ref{fig:cad_scheme}). These 
pressure values are measured with a catheter that is inserted in the coronary arteries;
however, to guarantee a sufficiently high accuracy, the patient needs to be put in pharmacologically-induced stress conditions.
This makes the procedure quite invasive and expensive. Other indices that can be measured under stress-free conditions 
have been proposed, like, for instance the \textit{instantaneous wave-free ratio} (iFR) \cite{gotberg2017instantaneous}.
For simplicity, we refer here to the FFR, but the same considerations apply to other indices.

Interventional cardiologists identify from images the location of a potential occlusion in the coronaries and 
measure with special catheters the pressure drop across it to quantify the FFR according to  equation \eqref{eq:ffr}. 
If the FFR is $< 0.8$, the patient is evaluated for a percutaneous intervention (PCI) or more invasive surgery, like a bypass graft.
PCI consists in the deployment of a wireframe structure called \textit{stent} along the occlusion
to restore the artery's physiological cross section.  To avoid the burden of measuring the FFR, 
HeartFlow \cite{HF}
launched a product to estimate the FFR with 
CFD. Whether measured or computed, the FFR must be accurate
for the well-being of the patient and for an optimal use of the resources, 
since false positives lead to unnecessary procedures and, more importantly, false negatives may lead to adverse events.

This assessment of CAD by FFR is designed for an effective short-term horizon, 
to decide treatment in the presence of pain. 
However, CAD is intrinsically multiphysics, combining fluid dynamics and tissue biology, and multiscale in time, as the effects of the blood flow occurring on the scale
of seconds develop CAD over months or years.
Biomedical engineers and clinicians have searched for
quantities with a long-term predictive value. 
The general consensus is that 
the WSS is a strong indicator of potential CAD.
For instance, in \cite{kumar2018low} it is demonstrated that low WSS in patients with non-obstructive CAD 
(high FFR)
is associated with severe Endothelial Dysfunction. In other words, low WSS ($< 1$ Pa) is associated with a marker of CAD, often 
anticipating plaque formation. On the other hand, \cite{kumar2018high_7873} 
shows that in patients with low FFR, so with a significant stenosis, high WSS ($> 4.71$ Pa) upstream the occlusion 
is associated with myocardial infarction. The analysis of the proximal WSS enhances the predictive value of the FFR.
In \cite{TUFARO2025119099}, the authors show that 
for patient having FFR around 0.8 (when the uncertainty can easily lead to a wrong assessment of the patient), adding a fast computation of the WSS can significantly 
improve the mid- to long-term risk stratification.

At this point, the following critical question arises: how to estimate reliably and rapidly the WSS?
When using CFD on patient-specific geometries, 
the WSS can be computed from the velocity field $\bu$. Optimized CFD
methods on specific workstations, starting from {Digital Subtraction Angiography} can provide a solution in around 10 minutes,
like in \cite{TUFARO2025119099}. On the other hand, CFD can be used to inform neural networks for a rapid assessment of the WSS, 
like in \cite{lin2025towards}.
In the spirit of a sustainable DT with a possible frequent collection from (non)invasive devices,
here we propose a different approach, based on DA.

\section{A DT Framework to Enhance the Treatment of CAD}
\label{sec:dt}

In general, a DT framework in healthcare includes the following components:
\begin{enumerate}
    \item the \underline{patient}, who is the source of data and the target of decisions;
    \item the \underline{devices to collect data}, with different frequencies, accuracy, and invasivity;
    \item  the \underline{P2D}  (Physical-to-Digital) step 
to convert measures into clinically relevant, but not measurable, data;
\item the \underline{D2P} (Digital-to-Physical) step to project a patient's data in time to make predictions; 
\item the \underline{clinician} (or \textit{human-in-the-loop}), who receives the predictions and a set of possible actions, and makes informed decisions on the patient's care (e.g., scheduling visits, acquiring further measures, {and/or performing} surgery).
\end{enumerate}

The bidirectional feedback loop~\cite{nationalacademy} between the physical space (patient's diseased coronary arteries), and its digital representation, operating under uncertainty is at the core of the DT system. Fig.~\ref{fig:bio_twin_scheme} depicts the continuous flow of information between the patient and the digital space. The
reciprocal exchange of information occurs in the P2D and D2P steps:
\begin{itemize}
    \item[-] {\bf P2D:} In this direction, data collected during clinical visits are channeled into the DT. To achieve this, DA becomes pivotal. DA enhances the computational model's synchronization and accuracy. 
    \item[-] {\bf D2P:} In this direction, 
the predictions, or optimized control strategies, from the DT's simulations inform the treatment of a patient's CAD. By treatment, we mean a therapeutic plan and the scheduling of recurrent visits to monitor the evolution of the pathology in time.
The DT functions as a platform to predict outcomes and optimize interventions or clinical visits, during which further data are collected.
We envision a DT based on physics-based modeling, DA, Machine Learning (ML), numerical analysis, and Uncertainty Quantification (UQ).
\end{itemize}

\begin{figure}[htb!]
\centering
\begin{overpic}[width=.85\textwidth, grid=false]{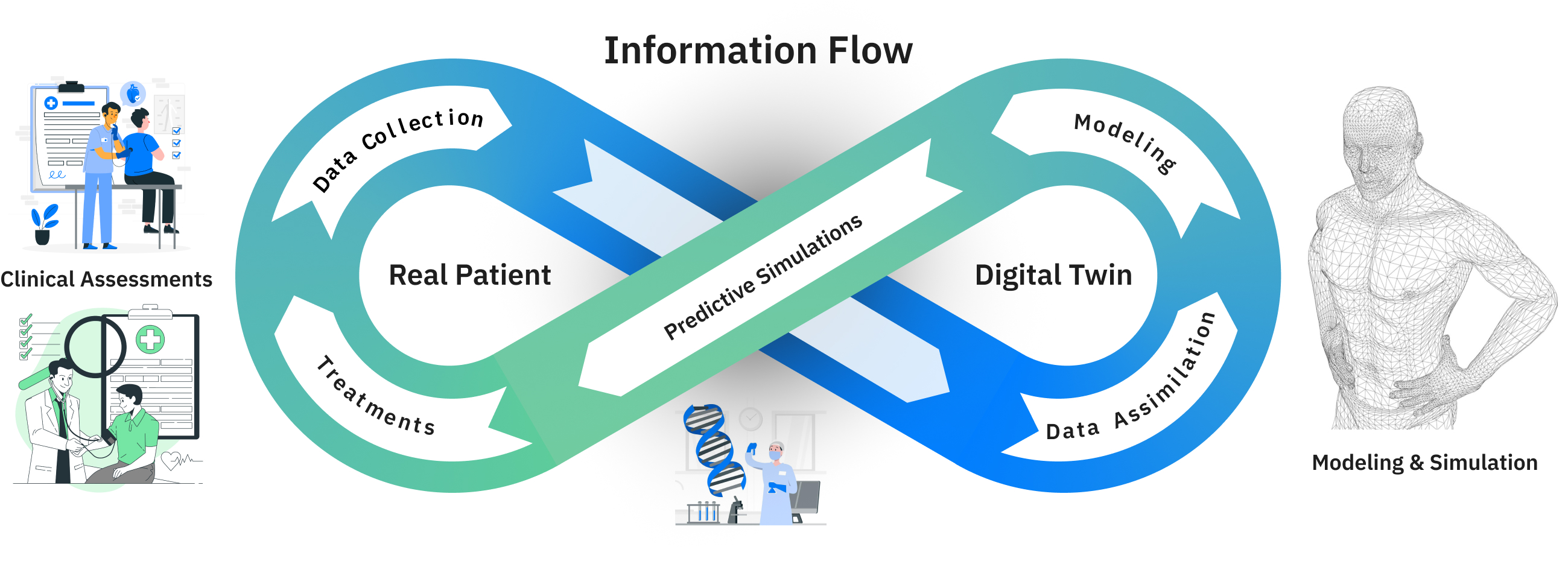}
\end{overpic}
    \caption{
    Information flow of a DT with bidirectional feedback loop.
    }
    \label{fig:bio_twin_scheme}
\end{figure}

To model the bidirectional coupling P2D-D2P, one could  adopt a PGM formulation successfully used for structural DTs of unmanned aerial vehicles~\cite{kapteyn2021probabilistic, tezzele2024adaptive} and civil engineering structures~\cite{torzoni2024digital}.
Several key sources of uncertainty are modeled in the PGM as time-dependent random variables: the underlying unknown physical state $S$, the observational data coming from regular clinical visits $O$, and the inferred digital state $D$. Additionally, a reward function $R$ helps to compute the optimal course of action $U$ for the duration of the treatment. The model determines the therapy according to the optimal policy that maximizes $R$. The edges within the PGM represent relationships between the connected variables via conditional probability distributions, which allows the uncertainty 
to be propagated throughout the entire graph, as displayed in Fig.~\ref{fig:pgm_complete}.

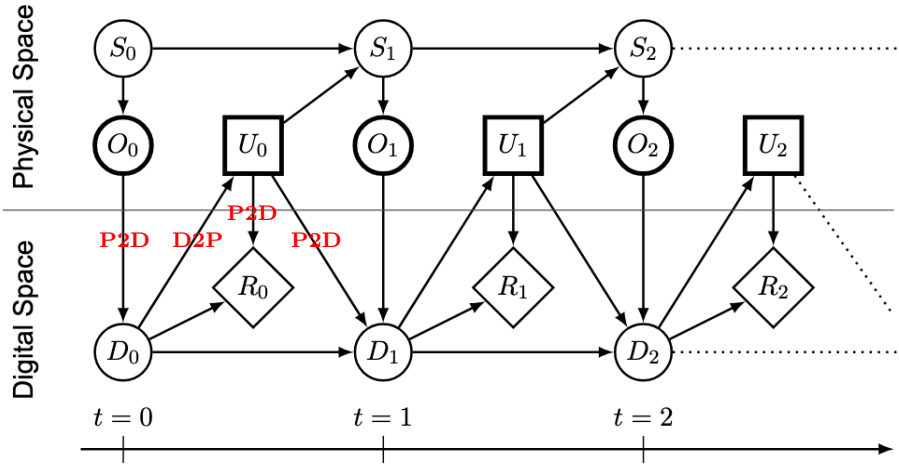
\begin{figure}[htb!]
\centering
\begin{overpic}[width=.8\textwidth, grid=false]{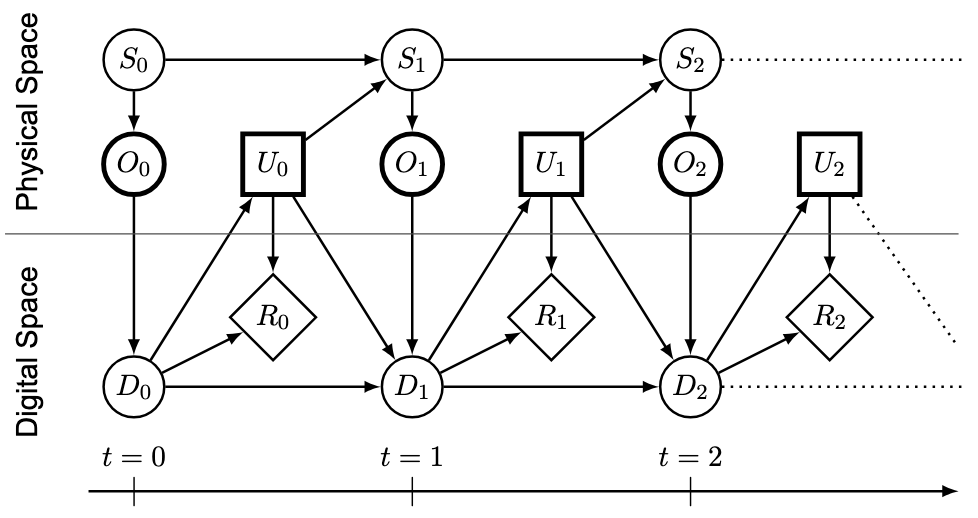}
\put(11,25){\textcolor{red}{\footnotesize {\bf P2D}}}
\put(19,25){\textcolor{red}{\footnotesize {\bf D2P}}}
\put(25,28){\textcolor{red}{\footnotesize {\bf P2D}}}
\put(32,25){\textcolor{red}{\footnotesize {\bf P2D}}}
\end{overpic}
    \caption{
    Example of a PGM with
    highlighted P2D and D2P couplings in the first time step.
Squares denote
actions, diamonds denote rewards, bold borders are for observed
(deterministic) quantities, and thin borders are for estimated (random) variables.
    }\label{fig:pgm_complete}
\end{figure}

For CAD, the goal is to estimate the hemodynamics indices ($S$ - the WSS in our case, see Sec.~\ref{sec:back}) 
from indirect observations $O$ (velocity) coming from US data. 
From these observations, we infer the hemodynamics indices through DA and physics-based numerical models 
(see Sec.~\ref{sec:VDA}), and we identify the digital state $D$. 
The ultimate goal is to minimize the risk of failure of the coronary, but also to optimize the number of future clinical visits, the prescriptions needed, or the overall cost of treatment. 
The clinicians represent the human-in-the-loop: they will select the optimal 
treatment among a list of possible actions, provided by the DTs framework {and} ranked based on the reward function.

In this general framework, we will focus herafter on the P2D and the D2P steps, which are the most critical 
for their mathematical aspects.

\section{P2D: Data Assimilation}
\label{sec:VDA}

The main target of our P2D step, at this preliminary stage of the construction of a DT system for CAD, is the accurate and rapid estimation of the WSS from possibly noninvasive velocity data. The velocity data collected by the US probes are generally sparse in a region of interest $\Omega$, and not necessarily on its boundary.
As mentioned earlier, numerical differentiation procedures for the data close to the boundary can achieve a WSS evaluation. Particular care must be taken with classical finite difference approaches, since the data are noisy and the differentiation can end up amplifying this noise. Pseudospectral techniques are capable of better noise management \cite{martin2024pseudo}.
Notice that the US probes do not necessarily retrieve all the velocity components, which introduces an additional error 
in the process. This data-driven approach is conceptually suboptimal, because it does not leverage our background knowledge of the physics, provided by the NSE 
(\ref{eq:pbc_nse})-(\ref{eq:pbc_nse-2}).
On the other hand,  the purely CFD-based approach can be used generally for retrospective clinical trials rather than in a DT framework.

We therefore propose to resort to DA for a reliable and sustainable personalization of the WSS evaluation. The basic idea,
illustrated in Fig. \ref{fig:da}, is to combine through DA the velocity data and the primitive variables of the NSE (velocity and pressure); then, the WSS can be retrieved by standard postprocessing procedures.

\begin{figure}[h]
    \centering
    \includegraphics[width=0.9\textwidth]{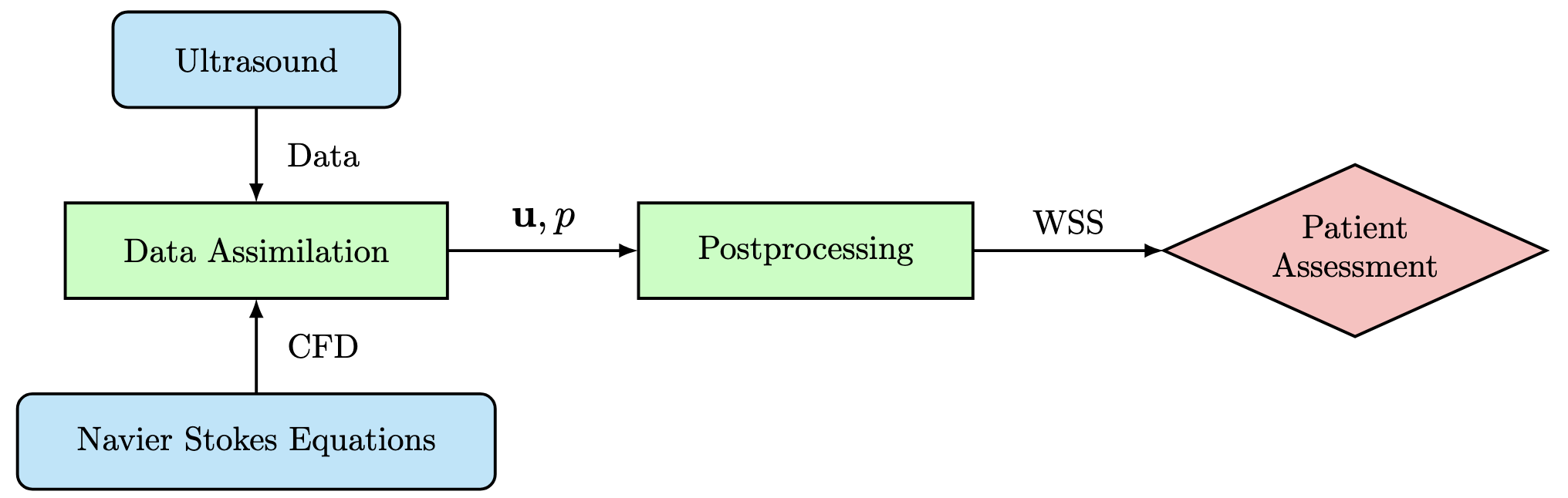}
    \caption{Block-diagram for the estimation of WSS through DA.}
    \label{fig:da}
\end{figure}

DA of the US data and the NSE can be clearly performed in different ways. One important criterion is the speed of execution: the DA has to be compatible with the clinical timelines and the DT framework, without reducing the accuracy of the estimation.

In the following, we consider three possible approaches.

\newcommand{\pin}{p_{\text{in}}}
\newcommand{\pout}{p_{\text{out}}}

\subsection{Variational Reduced-Order Modeling DA}
\label{sec:rom-da}

For the sake of simplicity, let us assume that the region of interest in the coronary at hand has only one inflow section and one outflow section, like in Fig. \ref{fig:coronary}. This assumption is realistic in many circumstances, but 
can be removed.
We also assume that the values $\pin$ and $\pout$ in \eqref{eq:pbc} and \eqref{eq:pbc-2}, respectively, are constant in space (commonly accepted in cardiovascular mathematics).

\begin{figure}
    \centering
    \includegraphics[width=0.5\linewidth]{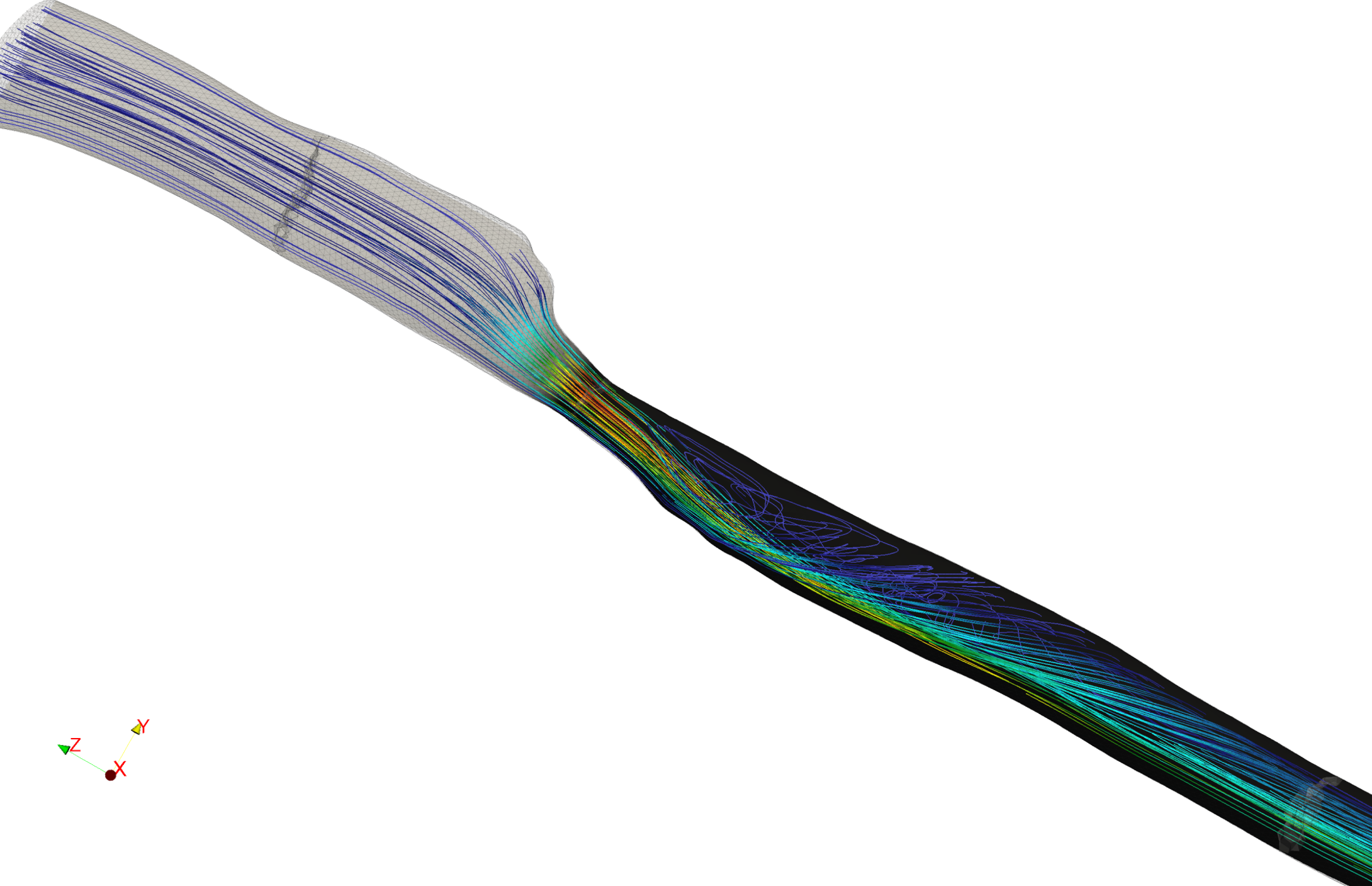}
    \caption{A coronary artery reconstructed from intravascular imaging. The plaque is evident at the end of the grey area. The region features one inlet (up, left) and one outlet (bottom, right).}
    \label{fig:coronary}
\end{figure}

One of the main challenges of cardiovascular simulations in a patient-specific setting is the (partial or entire) lack of information about the boundary conditions \cite{heywood1996artificial,veneziani98boundary,formaggia2002numerical,formaggia2010flow,quarteroni2016geometric,romarowski2018patient}. Consequently, the data $\pin$ and $\pout$ in \eqref{eq:pbc} and \eqref{eq:pbc-2} are seldom available in practice.
This is mainly due to the fact that measurement devices are invasive (and expensive if used routinely). 
On the other hand, available (velocity) data may be sparse in the region of interest $\Omega$.
A classical approach to perform DA in hemodynamics is therefore \textit{to consider the missing boundary conditions as the control variables to minimize the mismatch between the available data and the computed solution}. This {variational data assimilation (VDA)} approach was advocated years ago in a number of papers \cite{d2012variational,d2013uncertainty,bertoglio2012sequential,d2012applications,perego2011variational}. 
For instance, let $[\bd(\bx_l,t^m)]$, for $l\in[1,N_x]$, $m\in[1,N_t]$ be a set of velocity data over the time interval $[0,T]$ in $\Omega$. 
We consider the following functional: 
\begin{equation}\label{eq:vda}
\J (p_{\text{in}},p_{\text{out}}) \equiv \dfrac{1}{2} \sum\limits_{l,m}\| \bu (x_l,t^m; \pin,\pout) - \bd (x_l,t^m) \|^2, 
\end{equation}
where $\bu (x_l,t^m; \pin,\pout)$ is the solution of \eqref{eq:pbc_nse}, \eqref{eq:pbc_nse-2}, \eqref{eq:pbc}, and \eqref{eq:pbc-2}.
Then, a possible formulation of the DA problem reads: 
\begin{center}
find $\argmin \J(\pin,\pout)$,
\end{center}
which basically means that we minimize $\J(\bu(\pin,\pout),\bd)$ under the constraints
of \eqref{eq:pbc_nse}, \eqref{eq:pbc_nse-2}, \eqref{eq:pbc}, and \eqref{eq:pbc-2}.
This problem was considered in \cite{formaggia2008new} using a Lagrange multiplier approach. This means that the functional $\J$ is augmented by incorporating the constraints, weighted by the Lagrange multipliers. Eventually, we form the so-called Karush-Kuhn-Tucker (KKT) system by differentiating the augmented functional with respect to the state variables, the control variables, and the Lagrange multipliers (see, e.g., \cite{gunzburger_pfco}). 
The finite element discretization of this approach was analyzed in \cite{d2012variational}, while in \cite{d2013uncertainty} the Bayesian formulation of the solution
was introduced to quantify the
uncertainty in the results. We note that, while accurate, the methods in \cite{formaggia2008new,d2012variational,d2013uncertainty} are computationally too expensive to be used in clinical practice.

Other derivative-free approaches proceed by evaluating the functional $\J$ for different values of the control variables. This enables us to approximate $\J$ by fitting and
consequently to identify a guess of the optimal solution by minimizing the approximate 
functional (possibly after restricting the values that the control variables can take, by physical considerations). The procedure is iterated until a convergence criterion is fulfilled - see, e.g.,  \cite{cartis2019derivative,cartis2019improving,cartis2022escaping}.

In these approaches, each evaluation of $\J(\pin,\pout)$ requires the computation of the solution to \eqref{eq:pbc_nse}, \eqref{eq:pbc_nse-2}, \eqref{eq:pbc}, and \eqref{eq:pbc-2} for different values of $p_{in}, p_{out}$. This is easily be a computational bottleneck, since this may be needed for hundreds (or more) iterations. However, considering \eqref{eq:pbc_nse}, \eqref{eq:pbc_nse-2}, \eqref{eq:pbc}, and \eqref{eq:pbc-2} as a problem parametrized 
by the boundary data $\pin, \pout$, the computational cost can be drastically reduced by Reduced Order Models (ROMs).

ROMs (see, e.g., \cite{peter2021modelvol1, benner2020modelvol2, benner2020modelvol3, hesthaven2016certified, rozza2008reduced, quarteroni2011certified,quarteroni2015reduced} for reviews) 
have emerged as the methodology of choice to reduce the computational cost when Full Order Models (FOMs) have to be carried out for several parameter values. ROMs replace the FOM with a lower-dimensional approximation that captures the essential flow behavior. 
The reduction in computational time is achieved through a two-step procedure. 
In the first step, called the \textit{offline} stage, one constructs a database of several FOM solutions, also called snapshots, associated to given times and/or physical parameter values:
$\{ \bu_h^1, \ldots ,  \bu_h^M\}$, where $M$ is the total number of time instances or parameter values (or both).  
The FOM database is used to generate a ROM basis, $\{ \bphi_1, \ldots , \bphi_r \}$, where $r$ is the ROM dimension.
This ROM basis is (hopefully much) smaller than the high-dimensional FOM basis but still preserves the essential flow features.
A possible way to construct the reduced basis is Proper Orthogonal Decomposition (POD)~\cite{gunzburger_pfco,holmes2012turbulence,volkwein2013proper}, which stores the FOM solutions as columns of a matrix to which Singular Value Decomposition is applied. 
Then, the eigenvectors associated with the largest singular values (viewed as the most relevant information in the data) are extracted to generate a ROM basis. The effectiveness of the procedure is strongly related to how fast the singular values decay: the fastest 
the decay is, 
the fewer singular values are needed for the construction of the ROM, and more effective is the model reduction.

Next, the ROM basis is used in a Galerkin framework to construct the Galerkin ROM (G-ROM).
Specifically, we use the ROM velocity approximation $\bu_r(\bx,t) = \sum_{j=1}^{r} a_j(t) \bphi_j(\bx)$ to model the fluid velocity, $\bu(\bx,t)$, and project the resulting equations onto the ROM space, $\bX^r = \text{span} \{ \bphi_{1}, \ldots, \bphi_{r} \}$. 
The resulting G-ROM can be written as follows:
Given $\bur^{n}$, find $\bur^{n+1} \in \bXr$ such that $\forall \bvr \in \bXr$
\begin{eqnarray}
    \left( \frac{\bur^{n+1} - \bur^{n}}{\Delta t} , \bvr \right)
    + b^*(\bur^{n+1}, \bur^{n+1}, \bvr)
    + \nu (\nabla \bur^{n+1} , \nabla \bvr)
    = (\bff^{n+1},\bvr).
    \label{eqn:g-rom}
\end{eqnarray}
We note that the G-ROM~\eqref{eqn:g-rom} in this form does not include a ROM pressure approximation.
The reason is that we assume that the FOM snapshots are discretely divergence-free and, thus, so are the ROM basis functions.
As a result, the ROM approximation $\bur$ (which is a linear combination of basis functions) will also be discretely divergence-free, and the ROM pressure term in the weak form disappears. 
Alternative ROM formulations that include a pressure approximation can be found in, e.g., \cite{hesthaven2016certified,quarteroni2011certified,ballarin2015supremizer,decaria2020artificial}.

The model reduction is effective if the dimension of the G-ROM~\eqref{eqn:g-rom}, $r$, generally is orders of magnitude lower than the FOM dimension.

In the second step of the ROM strategy, called {online} stage, one uses the G-ROM~\eqref{eqn:g-rom} to quickly compute the solution for times and/or parameter values that are different from those used in the ROM construction.
The G-ROM~\eqref{eqn:g-rom} is an appealing alternative to the standard FOMs because its computational cost can be orders of magnitude lower than the FOM computational cost, while retaining a good level of accuracy.

ROMs have been used to reduce the computational cost in data assimilation~\cite{cao2007reduced,kaercher2018reduced,maday2015parameterized,stefanescu2015pod,popov2021multifidelity,donoghue2022multi}, and, specifically, in problems related to cardiovascular sciences, see, e.g., \cite{bertagna2014model,yang2017efficient,barone2020efficient}.  
It should be noted that, from a different perspective,  DA holds great promise in increasing ROM accuracy by incorporating available data into the ROM to account for both numerical and modeling errors~\cite{zerfas2019continuous,mou2023efficient}.

In particular, adding ROM for the 
rapid solution of \eqref{eq:pbc_nse}, \eqref{eq:pbc_nse-2}, \eqref{eq:pbc}, and \eqref{eq:pbc-2} for different values of the boundary data significantly accelerates 
the solution of the VDA assimilation problem by optimization. Thanks to the splitting between 
the offline and the online stage, preliminary (unpublished) results \cite{jimena2025} suggest that the optimization can be done in 15 minutes.
These results  demonstrate that the WSS based on DA in nontrivial geometries and 
with real US data affected by noise is more accurate than the 
WSS based only on data.

\subsection{Continuous Reduced-Order Modeling DA}
An alternative approach to VDA is {\it continuous DA (CDA)}, which is the insertion of coarse grain observational measurements directly into a mathematical model. CDA, which has connections with feedback control \cite{gunzburger_pfco} 
and 
nudging \cite{lakshmivarahan2013nudging}, has provided more accurate forecasts in applications ranging from the environmental to the medical sciences 
\cite{browning1998numerical,blomker2013accuracy,DAbook2016}.
To present this approach {in a clear, simplified setting}, consider the following general dynamical system, which can represent generically  a system of PDEs modeling 
cardiovascular flows:
\begin{eqnarray}
    \frac{\textrm{d} u}{\textrm{d} t}
    = F(u),
    \qquad
    u(0)
    = u_0,
    \label{eqn:cda-1}
\end{eqnarray}
where $u_0$ is the initial condition. 
We assume that data for the true solution $u$ are available on a coarse mesh, of size $H$, which is significantly larger than the fine mesh, of size $h$, on which the spatial discretization of eq.~\eqref{eqn:cda-1} is effected.
Let $I_H(u)$ be an available interpolation operator on the coarse mesh, or the orthogonal projection onto Fourier modes of wavenumbers smaller than $1/H$.  
The CDA strategy proposed in \cite{azouani2014continuous} is as follows:
\begin{eqnarray}
    \frac{\textrm{d} v}{\textrm{d} t}
    = F(v)
    + \mu (I_H(u) - I_H(v)),
    \qquad
    v(0)
    = v_0,
    \label{eqn:cda-2}
\end{eqnarray}
where $\mu$ is a nudging parameter. When the dynamical system in Eq.~\eqref{eqn:cda-1} is represented by Eqs.~\eqref{eq:pbc_nse}--\eqref{eq:pbc_nse-2} in 2D, it was proven in \cite{azouani2014continuous} that, {\it for any initial data $v_0$}, 
   $ \| v(t) - u(t) \|_{L^2}
    \longrightarrow 0$
{\it exponentially in time}.
Thus, even without knowing the initial data, $u_0$, the solution $u$ can be approximately reconstructed for large times. This is
one of the main advantages of CDA. 
Unfortunately, also in the case of CDA,
the computational cost limits its use in real clinical cases. 

Again, a possible workaround is represented by ROM. 
To date, the only CDA application at the ROM level was proposed in \cite{zerfas2019continuous}, where the CDA-ROM was applied to 2D flow past a circular cylinder at Reynolds numbers $Re=500$ and $Re=1000$ (a range that is significant for CAD in the presence of occlusions).
We note that the computational setting used in the numerical investigation in \cite{zerfas2019continuous} poses several challenges to the standard G-ROM~\eqref{eqn:g-rom}:
First, the G-ROM is used in the {\it under-resolved} regime, i.e., the number of ROM basis functions utilized to construct the G-ROM (e.g., $r=8$) is not enough to accurately represent the underlying complex dynamics.
Second, the  G-ROM is used in the {\it convection-dominated} regime, which poses significant challenges to reduced-order modeling.
The third challenge is that, in our numerical investigation, we built the G-ROM by using inaccurate snapshots.
Specifically, in \cite[Section 4.3]{zerfas2019continuous}, only $64\%$ of one period of FOM data to construct the ROM basis was used.
We emphasize that the under-resolved, convection-dominated, and data-scarce regimes are characteristic to the numerical simulation of complex flows in engineering, scientific, and medical applications.  
All of these above challenges can be tackled with a CDA-ROM.
In \cite{zerfas2019continuous},
the FOM snapshots are generated by using a finite element (FE) discretization with $(P_2, P_1^{\text{disc}})$ Scott–Vogelius elements on a barycenter refined Delaunay mesh with $103{,}000$ velocity degrees of freedom, and a BDF2 time discretization with a time step of $\Delta t = 0.002$. 
We generate a (coarse) mesh with size $H \gg h$ using the points $\{ \bx_l \}_{l = 0}^{N_x}$ where the velocity data $[\bd(\bx_l,t^m)]$, with $l\in[1,N_x]$, $m\in[1,N_t]$, are collected. 
To construct the CDA-ROM, we add a nudging term to the 
G-ROM~\eqref{eqn:g-rom}: 
Given $\bur^{n}$, find $\bur^{n+1} \in \bXr$ such that, $\forall \; \bvr \in \bXr$, we have: 
\begin{eqnarray}
    \left( \frac{\bur^{n+1} - \bur^{n}}{\Delta t} , \bvr \right)
    + b^*(\bur^{n+1}, \bur^{n+1}, \bvr)
    + \nu (\nabla \bur^{n+1} , \nabla \bvr)
    + \\ \mu \left( I_H (\bur^{n+1}) 
    - \bd(\cdot,t^{n+1}) , I_H(\bvr) \right)
    = (\bff^{n+1},\bvr), 
    \label{eqn:cda-rom}
\end{eqnarray}
where $I_H$ is an interpolation operator on the coarse mesh, and $b^*$ is the explicitly skew-symmetrized trilinear form.
In \eqref{eqn:cda-rom}, with a slight abuse of notation, $\bd(\cdot,t^{n+1})$ denotes a function that interpolates the data over the mesh of size $H$ at time $t^{n+1}$. 
The role of the nudging term is to steer the ROM velocity, $\bur$, towards the true solution, which corresponds to the velocity data, $[\bd(\bx_l,t^m)]$. 

The CDA-ROM \eqref{eqn:cda-rom} yields more accurate results than the G-ROM~\eqref{eqn:g-rom}
in the numerical simulation of a 2D flow past a cylinder at Reynolds numbers $Re=500$ and $Re=1000$ \cite{zerfas2019continuous}. This range of Reynolds numbers is relevant to coronaries.
Specifically, the cases $r=8$ and $r=16$ were considered, and compared the two ROMs with respect to several criteria: error, lift and drag coefficients, and kinetic energy.
The CDA-ROM~\eqref{eqn:cda-rom} is significantly more accurate than the standard G-ROM~\eqref{eqn:g-rom} with respect to all the criteria considered and for both Reynolds numbers.   
For $Re=500$, the case of snapshots generated over less than one period (i.e., $64\%$ and $84\%$ of one period) 
showed that DA can significantly 
improve the ROM accuracy when insufficient data is available to build the ROM (see, e.g., \cite[Fig.7]{zerfas2019continuous}). Furthermore, the improvement in the DA-ROM accuracy over the standard G-ROM accuracy was larger when inaccurate snapshots were used instead of accurate snapshots. 
In \cite[Section 4.4]{zerfas2019continuous}, an {\it adaptive} CDA-ROM, in which the nudging parameter is adjusted so that the CDA-ROM energy matches the true solution energy, was proposed too. 
The numerical results showed that the adaptive nudging slowly reduced the error and yielded a very accurate adaptive DA-ROM energy.

The error analysis of this CDA-ROM~\eqref{eqn:cda-rom} led to the following stability result
\cite[Lemma 3.3]{zerfas2019continuous}:
\begin{eqnarray}
    \| \bu_r^{n} \|
    \leq C_{\text{data}},
    \label{eqn:cda-rom-lemma}
\end{eqnarray}
where $C_{\text{data}}$ is a constant that depends on the input data. 
Furthermore \cite[Theorem 3.5]{zerfas2019continuous}, for an appropriately chosen interpolation operator $I_H$, it is possible to prove an {\it a priori} bound of the form
\begin{eqnarray}
    \| \bu^{n+1} - \bu_r^{n+1} \|^2
    \leq \ C^{n+1} \ \| \bu^{0} - \bu_r^{0} \|^2 
    + \text{ FOM}
    + \text{ ROM},
    \label{eqn:cda-rom-theorem}
\end{eqnarray}
where FOM (ROM) stands for the discretization error of the full (reduced) model,  $\bu^{0}$ and $\bu_r^{0}$ are the initial velocities, and $C$ is a constant. 
Since $|C| < 1$, the first term in the error bound \eqref{eqn:cda-rom-theorem} goes to zero {\it exponentially fast in time}. 

All these results sugges that CDA-ROM can be an effective tool in the design of our DT system.
However, extending the CDA-ROM \eqref{eqn:cda-rom} from the simplified setting
considered in \cite{zerfas2019continuous}
to the complex CAD setting poses several challenges. 
These include complex 3D geometries that need to be parametrized, accuracy enhancement, scarce available data, and convection-dominated flows. 
Strategies for tackling these challenges are outlined in Section~\ref{sec:road}.

\subsection{Physics-Informed Neural Networks in DA}
\label{sec:DA-PINN}

Physics-Informed Neural Networks (PINNs) are a declination of scientific ML based on the construction of a Neural Network (NN) that minimizes the residual of a PDE of interest. 
See, e.g., \cite{cai2021physics} and references therein.
If we generically denote problem \eqref{eq:pbc_nse}--\eqref{eq:pbc_nse-2} written in a residual form by $\mathcal{R}(\bu,p)=0$ and, similarly, the boundary conditions by $\B(\bu,p)=0$, the PINN solution of the flow problem of interest can be written as: find the parameters
$\boldsymbol{\vartheta}$ of a NN with input $(\mathbf{x},t)$ and output $(\bu(\mathbf{x},t), p(\mathbf{x},t))$ to minimize the functional
\begin{equation}\label{eq:PINN}
\begin{array}{l}
\J_{\text{PINN}} (\boldsymbol{\vartheta}) = w_1 \sum_{j,k} | \mathcal{R}(\bu(\mathbf{x}_j,t_k;\boldsymbol{\vartheta}), p(\mathbf{x}_j,t_k;\boldsymbol{\vartheta}))|^2 + \\[5pt]
\qquad \qquad \qquad w_2 |\B(\bu(\mathbf{x}_j,t_k;\boldsymbol{\vartheta}), \; p(\mathbf{x}_j,t_k;\boldsymbol{\vartheta})))|^2, 
\end{array}
\end{equation}
where $w_1$ and $w_2$ are positive weights, and the sum is extended to an appropriate number of points in space and time.
This approach may lead to fast solutions, but the accuracy depends on the design of the NN in a way that is not yet completely understood. 
Approximation results are available \cite{cybenko1989approximation,hornik1991approximation,guhring2020expressivity}, but their practical impact has yet to be fully exploited.
Therefore, the competitiveness of PINN over traditional CFD methods, e.g., finite element methods, 
is under scrutiny (see, e.g., \cite{grossmann2024can, rezaei2022mixed}). 
However, PINNs  provide a natural framework for DA
as one can incorporate the data mismatch in the solution of the problem by replacing the functional in \eqref{eq:PINN} with:
\begin{equation}\label{eq:PINNDA}
\J_{\text{PINN}}^{\text{DA}} (\boldsymbol{\vartheta})= \J_{\text{PINN}} (\boldsymbol{\vartheta})  + w_3  \sum\limits_{l,m}\| \bu (x_l,t^m; \boldsymbol{\vartheta}) - \bd (x_l,t^m) \|^2 .
\end{equation}
Note that the points where we minimize the residual and the points where we have measures 
are generally unrelated.

Preliminary (unpublished) results \cite{jimena2025} indicate that the DA-PINN significantly improves the accuracy of a standard PINN approach.
For instance, for the classical 2D Kim and Moin analytical solution \cite{kim1985application}, the errors (in the $L^2$ norm) without and with DA 
are displayed in Tab. \ref{tab:pinn}. The results (obtained with the python library DeepXDE \cite{lu2021deepxde}, using the same neural network design in the two cases) clearly indicate that the addition of data to the residual functional greatly improves the results.
It is also interesting to note that the pressure estimation greatly benefits from the data, which however refer to the velocity only.
Again, the WSS computation after the NSE is solved is just a standard postprocessing step based on numerical differentiation.
These results show that DA-PINN may be a viable option for the rapid and accurate estimation of the WSS, but 
deeper investigations on more complex and patient-specific geometries are in order.

\begin{table}
    \centering
    \begin{tabular}{c||c|c}
         &  \ $L^2$ errors PINN \ & \ $L^2$ errors DA-PINN \\
    \hline
    $u_1$ & 0.074 & 0.0023 \\
         $u_2$ & 0.13 & 0.0026 \\
         $p$ & 0.47 & 0.0022 \\
    \end{tabular}
    \caption{Comparison of the errors for PINN and DA-PINN for the Kim and Moin 2D problem. }
    \label{tab:pinn}
\end{table}

\begin{remark}
Aligning with the development of a DT system, we also aim to incorporate patients who have undergone PCI and
the deployment of a stent.  
This is a nontrivial task \cite{lefieux2022semi} - see Fig. \ref{fig:stent}.
Importantly, the struts of a stent can notably affect the WSS on endothelial tissue still in contact with blood flow \cite{ladisa_review_2022}. Specifically, a recent CFD study comparing 30 patient-specific geometries \cite{shah2025impact} revealed that WSS patterns result from a combination of macroscale aspects, such as vessel curvature and torsion, and microscale factors like the stent's shape and size. When curvature or torsion are minimal, the stent's influence on WSS is heightened. 

\begin{figure}
    \centering
    \includegraphics[width=0.48\linewidth]{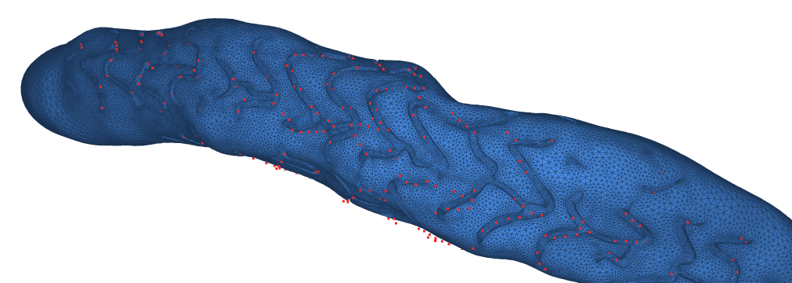} \hfill
        \includegraphics[width=0.48\linewidth]{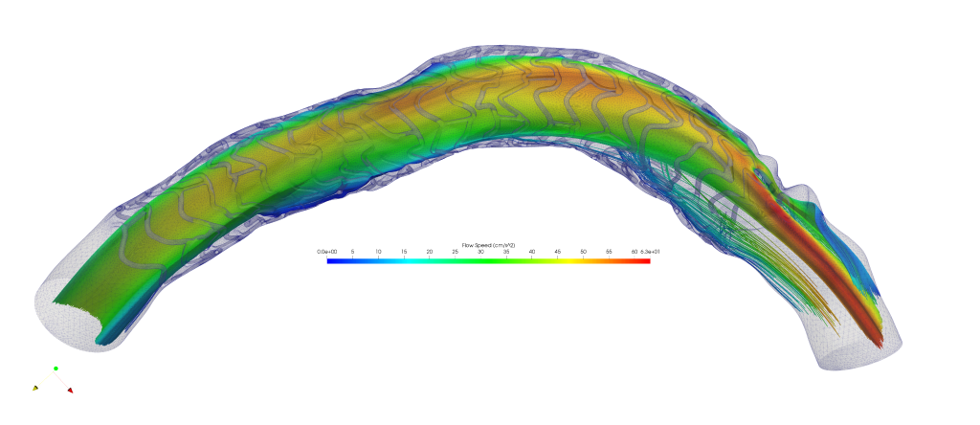}
    \caption{Lumen of a patient-specific coronary artery with the footprint of the stent (left) and blood flow streamlines (right).}
    \label{fig:stent}
\end{figure}

This underscores the necessity for precise and semi-automated methods in reconstructing patient-specific stented geometries, applicable to CFD (using advanced polytopal methods, for example \cite{antonietti2022virtual}) and, consequently, DA. Addressing these factors will present a formidable challenge in computational mathematics in the forthcoming years.   
\end{remark}

\section{D2P: From Data to Decision-Making}
\label{sec:PGM}

Scalable predictive DTs were recently introduced in \cite{kapteyn2021probabilistic} representing a DT with a graph. 
Graphical models \cite{jordan1999learning} combine probability theory with graph structures to help address two key challenges common in engineering and applied sciences: uncertainty and complexity. They tackle these challenges by breaking down complex systems into smaller, manageable components that are then mathematically connected in a consistent way. This modular setup makes it easier to integrate data, ensure internal consistency, and design scalable algorithms.

Fig.~\ref{fig:pgm_complete} illustrates such a graphical structure for a DT. 
Each time step can include new observational data and updated decisions, such as changes in control inputs. 
These updates are informed by models, which in our case vary in resolution and accuracy. 
A major challenge is that effective decision-making in DTs requires not just knowing the current system state, but predicting how the system will behave in the future --- ideally with quantified confidence. Embedding full physical models (e.g., PDEs) into the graph structure can support this, but at the cost of high computational expense, especially for complex  systems like those in our simulations. Hence, the importance of techniques aimed at reducing the computational cost is clear, as illustrated earlier. 
The proposed PGM allows us to monitor the patient, assimilate new observations, estimate the risk indices based on the personalized hemodynamics, and prescribe optimal treatments by comprising both physics-based models and data with quantified uncertainty. Moreover, the part of the graph modeling the digital space can incorporate future state predictions of the patient state, thus providing critical insights to clinicians, who make the actual decisions for the CAD treatment.

\subsection{PGM with Memory}

The main assumption of existing works using PGMs for DTs~\cite{kapteyn2021probabilistic, torzoni2024digital, tezzele2024adaptive, henaogarcia2025scitech, cotoarbua2024probabilistic}
is the Markovianity of the random variables, 
which means that the variables at time $t+1$ depend only on the variables at the previous time $t$. For example, following the arrows in Fig.~\ref{fig:pgm_complete}, the probability associated with $D_{t+1}$ depends only on $D_{t}$, on the realization of $U_{t}$, denoted by $u_t$, and on the observations $O_t = o_t$ as follows: $p(D_{t+1} \mid  D_{t}, U_{t} = u_{t}, O_{t+1} = o_{t+1})$.
This strong assumption can significantly decrease the accuracy of predictive digital models used for decision-making. 
If we consider a narrowed coronary with a growing plaque, it is reasonable to expect that the rate of change in the plaque size is as important as the plaque size itself. An abrupt change in a relatively low-risk patient would result in immediate treatment, while steady conditions in a compromised coronary artery could be treated with less aggressive solutions.
Therefore, to enable the use of PGM in clinical settings, it is important to weaken the Markovianity assumption and allow multiple variables from several past time instants to influence the predictions for the future digital state (Multi-Step Markov Chains)
and, hence, the corresponding treatment. 
By changing the topology of the graph with respect to the one in Fig.~\ref{fig:pgm_complete}, we expect improved DT's predictive accuracy which leads to better decision-making. In this new framework, the probability associated with $D_{t+1}$ will depend on the previous $k$ digital states and treatments, and on the last observation as follows: 
\begin{equation*}
p(D_{t+1} \mid  D_{t}, D_{t-1}, \dots, D_{t-k+1}, U_{t} = u_{t}, U_{t-1} = u_{t-1}, \dots, U_{t-k+1} = u_{t-k+1}, O_{t+1} = o_{t+1}).
\end{equation*}
At the clinical level, involving more steps to evaluate the rates of change in the process is very innovative. A possible drawback could be a degradation in computational efficiency
and scaling due to the need to store more data and the resulting factorization. 
Any potential scaling issue could be overcome by leveraging different information sources, such as data being assimilated with different fidelity levels or at different time scales,
as it is the case in our application. As we noted, data will eventually be collected from wearable and portable devices (with high frequency and, generally, low fidelity)
as well as from noninvasive and invasive measurements (lower frequency and higher fidelity).
This represents a critical direction for future research.

Another aspect to consider is that we often have limited clinically relevant data, especially compared to the number of variables we need to estimate. This imbalance makes inference and control particularly difficult using conventional methods. These obstacles suggest that we need new, scalable DT frameworks that embed physical models in a way that is both accurate and computationally efficient. We need methods that feature a rigorous treatment of uncertainty and can treat complex geometries or interacting physical phenomena \cite{gao2022physics,vadeboncoeur2023fully}. 
A robust DT system must jointly address: (i) Data assimilation updating system models using available measurements, and 
(ii) Optimal control making informed decisions to steer the system toward a desired goal.
These tasks are deeply coupled and must be addressed together in any practical DT application.

\subsection{Multifidelity PGM and UQ}

The new DT framework outlined above can be viewed as a hierarchical graphical modeling framework that integrates information from multiple sources of varying fidelity, accuracy, and computational cost. This framework is structured as a directed graph in which nodes represent submodels of different fidelities, and the outputs of some nodes serve as inputs to others. We embed both existing solvers and data sources, enabling a principled combination of models and data. Rather than discarding established models, we combine them at levels where they are most effective. 
We note that the traditional multifidelity approach, which uses lower-fidelity models to correct approximations in higher-fidelity ones, assumes that information sources can be hierarchically ordered based on predictive accuracy. However, in the presence of noisy data and model-form uncertainties, this assumption may not hold. We therefore envision a more flexible approach that learns nonlinear manifolds, which can be viewed as curved coordinate systems that encapsulate each fidelity separately. These nonlinear manifolds enable fast and accurate predictions even with limited simulation or sparse experimental data.

In the context of single-fidelity modeling, the goal is typically to represent a physical system using available data or imperfect models. Computational models approximate the desired solution by parameterizing an ansatz, which may take the form of classical mesh-based finite volume or finite element functions, or modern data-driven approaches such as deep neural networks \cite{vadeboncoeur2023fully}. Given limited data and computational resources, these approximations are inherently imperfect. To quantify their accuracy, appropriate error or risk metrics, such as mean squared error or statistical risk functions, must be minimized, typically through optimization over uncertain parameters and sample averages. However, challenges remain: optimal architecture selection, determination of collocation point density, and high computational costs associated with training. 

A critical gap in the literature involves handling operator uncertainty arising from parameterization. Most current approaches conduct operator learning for fixed parameterizations and treat the resulting operator as static, neglecting uncertainty in parameterization itself. This highlights the need for new methods that explicitly incorporate parameterization uncertainty from the outset. Existing single-fidelity methods often operate under deterministic settings where parameters are fixed; see, for instance, deep neural networks trained to minimize error at specific data points, as in PINNs \cite{raissi2019physics, karniadakis2021physics, cuomo2022scientific}. Operator learning techniques such as Fourier neural operators \cite{li2020fourier, wen2022u} and DeepONets \cite{lu2021learning, kovachki2023neural} aim to predict domain-wide solutions given initial or boundary conditions. Other approaches include neural ordinary differential equations \cite{chen2018neural, lee2021parameterized} and probabilistic methods for forward-inverse PDE mapping \cite{vadeboncoeur2023fully}. While effective, these methods typically do not address the complexities introduced by multifidelity systems.

To broaden these capabilities, we adopt a comprehensive definition of multifidelity that spans a wide range of model types, e.g., from high-fidelity PDEs to projection-based ordinary differential equations, implemented using deep neural architectures, Gaussian processes, and other statistical surrogates. Our framework encompasses varying levels of spatial and temporal resolution and accuracy, extending to multiscale and multiphysics systems. Here a modular multifidelity strategy can be formulated to overcome the computational burdens of high-fidelity-only models. By fusing data and model outputs from a heterogeneous ensemble of sources, we aim to construct a unified representation that improves prediction accuracy, especially in scenarios with sparse simulation or experimental data. The value of this multifidelity perspective, particularly in the absence of standardized validation datasets was recently emphasized in \cite{ren2023superbench}.  

Such a multifidelity framework can be built on advances in uncertainty quantification, particularly the use of high-fidelity data to correct local low-fidelity approximations --- an idea that has reduced computational expense in trust-region-based optimization \cite{alexandrov1998trust}. Discrepancy modeling has evolved to support more sophisticated corrections \cite{peherstorfer2018survey}, grounded in the premise that models can be hierarchically ranked by predictive utility \cite{kennedy2000predicting}. Several studies have developed surrogate models that exploit such hierarchies to integrate increasing levels of physics and numerical refinement \cite{forrester2007multi, narayan2014stochastic, le2014recursive, perdikaris2015multi, zhou2020generalized, chen2022multi}. However, this assumption often fails when a clear cost-benefit ranking among models is unavailable. Alternative approaches \cite{zhang2021multi, liu2019multi, de2020transfer, gorodetsky2020mfnets, meng2020composite, meng2021multi, pawar2022towards} have tried to relax these assumptions, but suffer from limitations such as fidelity inter-dependencies, constrained parameterizations, and the need for common variables across fidelities.

To overcome these limitations, we envision to utilize a multifidelity graph network architecture consisting of multiple elementary sets, each containing submodels derived from either physics or data. These are organized as directed graphs where nodes correspond to fidelity-specific models, and child nodes (higher fidelities) depend on parent nodes (lower fidelities). This structure mitigates the bottlenecks common in traditional multifidelity models. Such framework also incorporates empirical and reduced-order models as legitimate sources of information, thus preserving domain knowledge that might otherwise be lost in purely data-driven methods. Furthermore, one can perform transitioning from simple to complex geometries using nonlinear embeddings, such as variational autoencoders, to facilitate the learning of geometric representations across fidelities. Rather than using the low-fidelity solution directly as input for a high-fidelity model, we learn a common latent space that unifies inputs and outputs across fidelity levels. Once projected into this space, information from one fidelity can be passed to another in a decoupled and flexible manner.

\section{Conclusions: an Applied Math Challenge}
\label{sec:road}

The construction of DTs in clinics will be an exciting and challenging prospect in the years to come.
The potential impact of this integration of data, mathematical models, and experience in the definition of personalized precision medicine programs cannot be underestimated. 
As noted in \cite{sel2024building,sel2025survey}, there are critical steps to be taken, ranging from personalization of the mathematical models to accurate validation and verification steps, accompanied by uncertainty quantification and numerical analysis that could certify the quality of the quantitative elaboration. Although the use of digital replicas of patient-specific apparatuses, thanks to the combination of data and models, already demonstrated great potential to understand the physiopathology of some diseases, the transition to ``twins'' has the potential to be a revolution in clinics, similar to the one brought about by the advent of imaging devices at the beginning of the 20$^{\text{th}}$ century.
From ``replicas'' to twins, from ``combination'' to ``integration'' of data and models: these steps 
challenge applied and computational mathematicians in the development of fast and reliable methods of personalization and synthesis. The decision-making  support demands for a reliable integration of data from different sources, collecting data with varying levels of reliability and frequency, and multi-physics, multi-scale (in time and space) models.\footnote{It is worth noting that improving the quality of data is also a challenge, involving more biomedical engineers and clinicians than mathematicians. Beyond the technological advancements of devices, we stress the need to create standard protocols for multi-center studies that could provide significant amounts of data with uniform collection procedures.}

DA is certainly a key step in this process, being central in the certified personalization of available models.
The advent of ROM methodologies enables the use of accurate techniques within clinical timelines.
However, there are still challenges that need to be overcome, like the construction of geometrical atlases to leverage, for instance, the offline stage performed for multiple patients. 
Another challenge is the development of accurate ROMs in convection-dominated regimes, which are relevant to CAD.
To tackle this challenge, a promising research avenue is the use of ROM spatial filtering stabilization (inspired by Large Eddy Simulation  
\cite{girfoglio2021pod,girfoglio2023hybrid,strazzullo2022consistency,tsai2025time,ivagnes2025data,quaini2024bridging}) in DA.

Yet another challenge is the setup of Scientific Machine Learning techniques like PINN for the identification and personalization of models, such that we gain understanding and confidence comparable to the ones we have now for traditional computational tools like Finite Elements.

Furthermore, the combination of different methods (for instance ROM and PINN to provide methodologies that leverage the best of the different approaches) should also be developed to make the personalization rapid and accurate.

Several challenges concern also the D2P step, particularly for the heterogeneity of data that will be eventually collected (in terms of frequency and fidelity).
Using PGM with memory (multistep chains) in the data elaboration seems to be a promising framework to be carefully validated and certified, with an accurate UQ analysis.  

Ideally, as noted in \cite{sel2024building}, the construction of a DT system should go beyond a single-disease object, and should consider each patient as a physiological system with an effective monitoring for the prevention of 
different types of pathologies. 
Strangely enough, in medical papers on this topic the role of mathematics is overlooked, in favor of a computer science/engineering focus. 
We believe that, on the contrary, the role of mathematics in this long-term (but hopefully not so long-term) objective is critical, and as much as the combination of general deductive reasoning and inferential knowledge required by DT is challenging, it is also exciting. 
The ultimate goal is the improvement of daily healthcare, where, we strongly believe, mathematicians can and must play a primary role at the frontline.

\section*{Acknowledgements}
AV, OS, and TI acknowledge the support of the US NSF under Award Numbers DMS-2012286, DMS-2345048, and DMS-2012253, respectively.
AV also acknowledges (i) the support of US NSF Award Number DMS2038118 (PI: J. Nagy),
(ii) the US NIH under Grant No. R01EB031101 (PI: B. Lindsey), 
(iii) the help of Jimena Martin-Tempestti (Emory University, USA) and Micol Bracco
(Politecnico di Torino, IT), and (iv) the inspiration provided by 
Brooks Lindsey (Georgia Institute of Technology, USA) and Habib Samady (Emory University, USA).

\bibliographystyle{abbrvurl}
\bibliography{the_one_and_only.bib}

\end{document}